\documentclass[11pt]{article}
\title{ON YAU RIGIDITY THEOREM FOR MINIMAL SUBMANIFOLDS IN SPHERES  \footnote{
\indent Research supported by the NSFC, Grant No.
11071211, 10771187; the Trans-Century Training Programme
\newline \indent Foundation for Talents by the Ministry of
Education of China.}}
\author{JUAN-RU GU AND HONG-WEI XU}
\date{}
\begin{document}

\maketitle
\begin{abstract}
In this note, we investigate the well-known Yau rigidity theorem for
minimal submanifolds in spheres. Using the parameter method of Yau
and the DDVV inequality verified by Lu, Ge and Tang, we prove that
if $M$ is an $n$-dimensional oriented compact minimal submanifold in
the unit sphere $S^{n+p}(1)$, and if
$K_{M}\geq\frac{sgn(p-1)p}{2(p+1)},$ then $M$ is either a totally
geodesic sphere, the standard immersion of the product of two
spheres, or the Veronese surface in $S^4(1)$. Here $sgn(\cdot)$ is
the standard sign function. We also extend the rigidity theorem
above to the case where $M$ is a compact submanifold with parallel
mean curvature in a space form.
\end{abstract}
 \section{Introduction}
\hspace*{5mm}It plays an important role in geometry of submanifolds to
investigate rigidity of minimal submanifolds. After the pioneering
rigidity theorem for closed minimal submanifolds in a sphere due to
Simons \cite{Simons}, a series of striking rigidity results for
minimal submanifolds were proved by several geometers \cite{Chern,
Lawson, Yau}. Let $M^n$ be an $n$-dimensional compact Riemannian
manifold isometrically immersed into an $(n+p)$-dimensional complete
and simply connected Riemannian manifold $F^{n+p}(c)$ with constant
curvature $c$. Denote by $K_M$ and $H$ the sectional curvature and
mean curvature of $M$ respectively. In 1975, Yau \cite{Yau} first
proved the
following celebrated rigidity theorem for minimal submanifolds in spheres under sectional curvature pinching condition.\\\\
\textbf{ Theorem A.} \emph{Let $M^n$ be an $n$-dimensional oriented
compact minimal submanifold in $S^{n+p}(1)$. If
$K_{M}\geq\frac{p-1}{2p-1},$ then either M is the totally geodesic
sphere, the standard immersion of the product of
two spheres, or the Veronese surface in $S^4(1)$. }\\\\
\hspace*{5mm}The pinching constant above is the best possible in the
case where $p=1$, or $n=2$ and $p=2$. It's better than the pinching
constant of Simons \cite{Simons} in the sense of the average of
sectional curvatures. Later, Itoh \cite{Itoh2} proved that if $M^n$
is an oriented compact minimal submanifold in $S^{n+p}(1)$ whose
sectional curvature satisfies $K_{M}\geq\frac{n}{2(n+1)},$ then $M$
is the totally geodesic sphere or the Veronese submanifold. Further
discussions in this direction have been carried out by many other
authors (see \cite{Ejiri, Kozlowski, Li, Shen, Xu, Xu1, Xu2}). An important problem is stated as follows.\\\\
\textbf{Open Problem B.} \emph{What is the best pinching constant
for the rigidity theorem for oriented compact minimal submanifolds
in a unit sphere under sectional (Ricci, scalar, resp.) curvature
pinching condition?}\\\\
\hspace*{5mm}Up to now, the problem above is still open. In
particular, Lu's conjecture \cite{Lu4}, a scalar curvature pinching
problem for minimal submanifolds in a unit sphere, has not been
verified yet. In this note, using Yau's parameter method \cite{Yau}
and the DDVV conjecture proved by Lu, Ge and Tang \cite{Ge, Lu2}, we
prove the following
rigidity theorem for minimal submanifolds in spheres.\\\\
\textbf{Theorem 1.} \emph{Let $M^n$ be an $n$-dimensional oriented
compact
 minimal submanifold in the unit sphere $S^{n+p}(1)$.
 If  $$K_{M}\geq\frac{sgn(p-1)p}{2(p+1)},$$
then M is either a totally geodesic sphere, the standard immersion
of the product of
two spheres, or the Veronese surface in $S^4(1)$. Here $sgn(\cdot)$ is the standard sign function.} \\\\
\textbf{Remark 1.} When $2<p<n$, our pinching constant in Theorem 1 is better than
ones given by Yau \cite{Yau} and Itoh \cite{Itoh2}.\\\\
\hspace*{5mm}More generally, we obtain the following rigidity result
for
submanifolds with parallel mean curvature in spaces forms.\\\\
\textbf{Theorem 2.} \emph{Let $M^n$ be an $n$-dimensional oriented
compact
 submanifold with parallel mean curvature $(H \neq 0 )$
in $F^{n+p}(c)$. If $c+H^2>0$ and $$K_M \ge
\frac{sgn(p-2)(p-1)}{2p}(c+H^2),$$
 then $M$ is either a totally
umbilical sphere $S^n(\frac{1}{\sqrt{c+H^2}})$ in $F^{n+p}(c)$, the
standard immersion of the product of
two spheres or the Veronese surface in $S^4(\frac{1}{\sqrt{c+H^2}})$.} \\\\
\section{Notation and lemmas}\hspace*{5mm}Throughout this paper,
let $M^{n}$ be an $n$-dimensional compact Riemannian manifold
isometrically immersed into an $(n+p)$-dimensional complete and
simply connected space form $F^{n+p}(c)$ of constant curvature $c$. We
shall make use of the following convention on the range of indices:
$$ 1\leq A,B,C,\ldots\leq n+p;\ 1\leq i,j,k,\ldots\leq n;\ n+1\leq
\alpha,\beta,\gamma,\ldots\leq n+p.$$ Choose a local field of
orthonormal frames \{$e_{A}$\} in $F^{n+p}(c)$ such that, restricted
to $M$, the $e_{i}$'s are tangent to \emph{M}. Let \{$\omega _{A}$\}
and \{$\omega _{AB}$\} be the dual frame field and the connection
1-forms of $F^{n+p}(c)$ respectively. Restricting these forms to
\emph{M}, we have
\begin{eqnarray}&&\omega_{\alpha i}=\sum_{j} h^{\alpha}_{ij}\omega_{j},  \, \,
h^{\alpha}_{ij}=h^{\alpha}_{ji},\nonumber\\
&&h=\sum_{\alpha,i,j}h^{\alpha}_{ij}\omega_{i}\otimes\omega_{j}\otimes
e_{\alpha},\,\,\xi=\frac{1}{n}\sum_{\alpha,i}h^{\alpha}_{ii}e_{\alpha},\nonumber
\\
&&R_{ijkl}=c(\delta_{ik}\delta_{jl}-\delta_{il}\delta_{jk})+\sum_{\alpha}(h^{\alpha}_{ik}h^{\alpha}_{jl}-h^{\alpha}_{il}h^{\alpha}_{jk}), \\
&&R_{\alpha\beta kl}=\sum_{i}(h^{\alpha}_{ik}h^{\beta}_{il}-h^{\alpha}_{il}h^{\beta}_{ik}),
\end{eqnarray}
where $h, \xi, R_{ijkl}, R_{\alpha\beta kl},$ and
$\overline{R}_{ABCD}$ are the second fundamental form, the mean
curvature vector, the curvature tensor, the normal curvature tensor
of $M$, and the curvature tensor of $N$, respectively. We define
$$S=|h|^{2}, \ H=|\xi|, \ H_{\alpha}=(h^{\alpha}_{ij})_{n\times
n}.$$
The scalar curvature $R$ of $M$ is given by
$$ R=n(n-1)c+n^{2}H^{2}-S.$$
Denote $K_{M}(p,\pi)$ the sectional curvature of $M$ for tangent
2-plane $\pi\subset T_{p}M$ at point $p \in M$. Set $K_{\min}(p) =
\min_{\pi\subset T_{p}M} K_M(p,\pi)$. From $\cite{Yau}$, we have the following lemma.\\\\
\textbf{Lemma 1.} \emph{If $M^n$ is a submanifold with parallel
mean curvature and positive sectional curvature in $F^{n+p}(c)$,
then M is a pseudo-umbilical submanifold. }\\\\
\hspace*{5mm}Let $M$ be a submanifold with parallel
mean curvature vector $\xi$. Choose $e_{n+1}$ such that it is parallel to $\xi$, and
\begin{equation}
trH_{n+1}= nH,  \ \ trH_\alpha= 0,\ \  \alpha \neq n+1.
\end{equation}
Set
\begin{equation}
 S_H=
trH^2_{n+1},\ \ S_I= \sum_{\alpha\neq n+1} trH^2_\alpha.
\end{equation}
 When $M$ is a pseudo-umbilical
submanifold, we have
\begin{equation}
 S_H= trH^2_{n+1}=nH^2.
 \end{equation} Denoting the first and second covariant derivatives of $ h^\alpha_{ij} $ by
$h^\alpha_{ijk} $ and $h^\alpha_{ijkl}$ respectively, we have
\begin{eqnarray*}
&&\sum_{k} h^\alpha _{ijk}\omega_k
 = dh^\alpha_{ij}-\sum_{k} h^\alpha_{ik}\omega_{kj}
 -\sum_{k} h^\alpha_{kj}\omega_{ki}-\sum_{\beta} h^\beta_{ij}\omega_{\beta\alpha},
\\
&&
 \sum_{l} h^\alpha _{ijkl}\omega_l= dh^\alpha_{ijk}-
 \sum_{l} h^\alpha_{ijl}\omega_{lk}-\sum_{l} h^\alpha_{ilk}\omega_{lj}-
 \sum_{l} h^\alpha_{ljk}\omega_{li}-\sum_{\beta} h^\beta_{ijk}\omega_{\beta\alpha}.
\end{eqnarray*}
Then we have
$$
h^\alpha_{ijk}=h^\alpha_{ikj},\ \  h^\alpha _{ijkl}-h^\alpha
_{ijlk}=\sum_{m} h^\alpha_{im}R_{mjkl}+
 \sum_{m} h^\alpha_{mj}R_{mikl}-\sum_{\beta} h^\beta_{ij}R_{\alpha\beta kl},
$$
\begin{equation}
 \Delta h^\alpha_{ij}=\sum_{k} h^\alpha _{ijkk}=\sum_{k} h^\alpha _{kkij}+
 \sum_{k}\Big(\sum_{m} h^\alpha_{km}R_{mijk}+\sum_{m} h^\alpha_{mi}R_{mkjk}-
 \sum_{\beta} h^\beta_{ki}R_{\alpha\beta jk}\Big).
 \end{equation}
 \hspace*{5mm}The following lemma will be used in the proof of our main results.\\\\
\textbf{Lemma 2(\cite{Yau}).} \emph{If $M^n$ is a submanifold with parallel mean curvature in $F^{n+p}(c)$, then
either $H\equiv0$ or H is non-zero constant and $H_{n+1}H_\alpha =H_\alpha
H_{n+1}$ for all $\alpha$.}\\\\
\hspace*{5mm}The DDVV inequality proved by Lu, Ge and Tang \cite{Ge,Lu2} is stated as follows.\\\\
\textbf{DDVV Inequality.} \emph{Let $B_{1},...,B_{m}$ be symmetric
$(n\times n)$-matrices, then
  \begin{equation}
\sum_{r,s=1}^m\|[B_{r},B_{s}]\|^2\leq \Big(\sum_{r=1}^m\|B_{r}\|^2\Big)^2,\end{equation}
where the equality holds if and only if under some rotation all $B_{r}$'s are zero except two
matrices which can be written as
$$\tilde{B}_{r}=P\left(\begin{array}{lllll}0&\mu&0&\cdots&0\\\mu&0&0&\cdots&0\\0&0&0&\cdots&0\\\vdots&
\vdots&\vdots&\ddots&\vdots\\0&0&0&\cdots&0\end{array}\right)P^t,\hspace*{10mm}
\tilde{B}_{s}=P\left(\begin{array}{lllll}\mu&0&0&\cdots&0\\0&-\mu&0&\cdots&0\\0&0&0&\cdots&0\\\vdots&
\vdots&\vdots&\ddots&\vdots\\0&0&0&\cdots&0\end{array}\right)P^t,$$ where P is an orthogonal $(n\times n)$-matrix.
Here $\|\cdot\|^ 2$ denotes the sum of squares of entries of the matrix and $[A,B]=AB-BA$ is
the commutator of the matrices A, B.}\\\\
\hspace*{5mm}For further discussions about the DDVV inequality, we refer to see \cite{DeSmet, Ge, Lu1, Lu2, Lu3, Lu4}.\\\\
\section{Proof of the theorems}
\hspace*{5mm}When $M^n$ be a minimal submanifold  in $S^{n+p}(1)$,
we have $trH_{\alpha}=0$ for all $\alpha$ and $\sum_{i}h^\alpha
_{iikl}= 0$. It follows from (6)  that
\begin{equation}
 \Delta h^\alpha_{ij}=\sum_{k,m }
h^\alpha_{km}R_{mijk}+\sum_{k,m } h^\alpha_{mi}R_{mkjk}-
\sum_{k,\beta} h^\beta_{ki}R_{\alpha\beta jk}.\end{equation} Thus
\begin{equation} \sum_{i,j,\alpha} h^\alpha _{ij} \Delta h^\alpha_{ij}
=\sum_{i,j,k,m,\alpha } (h^\alpha_{ij}h^\alpha_{km}R_{mijk}+
h^\alpha_{ij}h^\alpha_{mi}R_{mkjk}) - \sum_{i,j,k,
\alpha,\beta} h^\alpha_{ij}h^\beta_{ki}R_{\alpha\beta jk}.\end{equation}\\
\textbf{Proof of Theorem 1.} By using (1) and (2), we get
\begin{eqnarray*} &&\sum_{i,j,k,m,\alpha} h^\alpha_{ij}h^\alpha_{km}R_{mijk}+
 \sum_{i,j,k,m,\alpha}
 h^\alpha_{ij}h^\alpha_{mi}R_{mkjk} \\
    &=&nS +\sum_{\alpha,\beta} trH_\beta\cdot tr(H^2_\alpha H_\beta)-\sum_{\alpha,\beta}[tr(H_\alpha H_\beta)]^2
    -\sum_{\alpha,\beta} [tr(H^2_\alpha H^2_\beta)-tr(H_\alpha
    H_\beta)^2],\end{eqnarray*}
  and $$\sum_{i,j,k,
\alpha,\beta} h^\alpha_{ij}h^\beta_{ki}R_{\alpha\beta jk}=\sum_{\alpha,\beta} [tr(H^2_\alpha H^2_\beta)-tr(H_\alpha
    H_\beta)^2]. $$
    Since $(tr(H_\alpha H_\beta))$ is a symmetric
$(p\times p)$-matrix, we can choose the normal frame fields
$\{e_{\alpha}\}$ such that $$tr(H_\alpha H_\beta)= trH^2_\alpha
\cdot \delta_{\alpha\beta}.$$ This implies
\begin{equation} \sum_{\alpha,\beta} [tr(H_\alpha H_\beta)]^2=
  \ \sum_{\alpha}(trH^2_\alpha)^2.\end{equation}
  From above equalities, we obtain
    \begin{eqnarray}
    \sum_{i,j,\alpha} h^\alpha _{ij} \Delta h^\alpha_{ij}&
=&-anS+(1+a)\sum_{i,j,k,m,\alpha}(
h^\alpha_{ij}h^\alpha_{km}R_{mijk}+
h^\alpha_{ij}h^\alpha_{mi}R_{mkjk})\nonumber\\
    &&+(a-1)\sum_{\alpha,\beta} [tr(H^2_\alpha H^2_\beta)
    -tr(H_\alpha
    H_\beta)^2]+ a\sum_{\alpha,\beta} (trH_\alpha^2)^2,\end{eqnarray}
for all real number $a$.
For fixed $\alpha$, we choose the orthonormal frame fields $\{e_i\}$
such that $h^\alpha_{ij}=\lambda^\alpha_i\delta_{ij}$. Hence, we get
\begin{eqnarray}
&&\sum_{i,j,k,m}h^\alpha_{ij}h^\alpha_{km}R_{mijk}+
    \sum_{i,j,k,m}
    h^\alpha_{ij}h^\alpha_{mi}R_{mkjk}\nonumber \\
 &=&\sum_{i,k} \lambda^\alpha_i\lambda^\alpha_k R_{kiik}+\sum_{i,k} \lambda^\alpha_i\lambda^\alpha_i
 R_{ikik}\nonumber\\
&=&\frac{1}{2}\sum_{i,j} (\lambda^\alpha_i-\lambda^\alpha_j)^2R_{ijij}\nonumber\\
& \geq &\frac{1}{2}K_{\min} \sum_{i,j}(\lambda^\alpha_i-\lambda^\alpha_j)^2\nonumber\\
&=&nK_{\min} (trH^2_\alpha),\end{eqnarray}
which implies that
\begin{equation}\sum_{i,j,k,m,\alpha}h^\alpha_{ij}h^\alpha_{km}R_{mijk}+
    \sum_{i,j,k,m,\alpha}
    h^\alpha_{ij}h^\alpha_{mi}R_{mkjk} \geq nK_{\min}S.\end{equation}
  On the other hand, by a direct computation and the DDVV inequality, we obtain
  \begin{eqnarray}\sum_{\alpha,\beta}tr(H^2_\alpha H^2_\beta)
    -tr(H_\alpha
    H_\beta)^2&=&\frac{1}{2}\sum_{\alpha,\beta}tr(H_{\alpha} H_{\beta}-H_{\beta}
H_{\alpha})^{2}\nonumber\\
&\leq&\frac{1}{2}sgn(p-1)\Big(\sum_{\alpha}trH_{\alpha}^{2}\Big)^2\nonumber\\
&=&\frac{1}{2}sgn(p-1) S^2,\end{eqnarray}
 where $sgn(\cdot)$ is the standard sign function. It follows from (11), (13) and (14) that
  \begin{eqnarray}
   \frac{1}{2}\Delta S&=&\sum_{i,j,\alpha} (h^\alpha_{ijk})^2 +\sum_{i,j,\alpha} h^\alpha_{ij}\Delta
    h^\alpha_{ij}\nonumber\\
 &\geq&\sum_{i,j,k,\alpha}(h^\alpha _{ijk})^2 -anS+(1+a)nK_{\min}S
 +\Big[\frac{a}{p}+\frac{sgn(p-1)}{2}(a-1)\Big]S^2, \end{eqnarray}
for  $0\leq a<1$.
  Taking $a=sgn(p-1)\frac{p}{p+2}$, we get
$$
  \frac{1}{2}\Delta S \geq nS\Big[\Big(1+sgn(p-1)\frac{p}{p+2}\Big)K_{\min}-sgn(p-1)\frac{p}{p+2}\Big].$$
 It follows from the assumption and the maximum principal that $S$ is a
 constant, and
 $$S\Big[\Big(1+sgn(p-1)\frac{p}{p+2}\Big)K_{\min}-sgn(p-1)\frac{p}{p+2}\Big]=0.$$
If there is a point $q\in M$ such that
 $K_{\min}(q)>\frac{sgn(p-1)p}{2(p+1)}, $ then
$S=0$, i.e., $M$ is totally geodesic. If
$K_{\min}\equiv\frac{sgn(p-1)p}{2(p+1)}$, then inequalities in (13),
(14) and (15) become equalities. From the DDVV inequality we obtain
$p\leq 2$. This together with Theorem A implies $M$ is the product
of
two spheres or the Veronese surface in $S^4(1)$. This completes the proof of Theorem 1.\\\\
\hspace*{5mm} When $M^n$ is a  submanifold with parallel mean
curvature in $F^{n+p}(c)$, we have $\xi=He_{n+1}$, and
$\sum_{i}h^\alpha _{iikl}= 0$ for $\alpha\neq n+1$. It follows from
(6) and Lemma 2 that
\begin{equation}
 \Delta h^\alpha_{ij}=\sum_{k,m }
h^\alpha_{km}R_{mijk}+\sum_{k,m } h^\alpha_{mi}R_{mkjk}-
\sum_{k,\beta \neq n+1} h^\beta_{ki}R_{\alpha\beta jk},\ \ \ \ \
\alpha\neq n+1.\end{equation} Thus
\begin{eqnarray} \sum_{i,j,\alpha\neq n+1} h^\alpha _{ij} \Delta h^\alpha_{ij}
&=&\sum_{i,j,k,m,\alpha\neq n+1}
(h^\alpha_{ij}h^\alpha_{km}R_{mijk}+
h^\alpha_{ij}h^\alpha_{mi}R_{mkjk})\nonumber\\
&& - \sum_{i,j,k,
\alpha,\beta\neq n+1} h^\alpha_{ij}h^\beta_{ki}R_{\alpha\beta jk}.
\end{eqnarray}\\
\textbf{Proof of Theorem 2.} Applying (1) and (2), we get
\begin{eqnarray*} &&\sum_{i,j,k,m,\alpha\neq n+1} h^\alpha_{ij}h^\alpha_{km}R_{mijk}+
 \sum_{i,j,k,m,\alpha\neq n+1}
 h^\alpha_{ij}h^\alpha_{mi}R_{mkjk} \\
    &=&ncS_I +\sum_{\alpha \neq n+1,\beta} trH_\beta\cdot tr(H^2_\alpha H_\beta)-\sum_{\alpha\neq
   n+1,\beta}[tr(H_\alpha H_\beta)]^2 \\
    &&-\sum_{\alpha,\beta\neq n+1} [tr(H^2_\alpha H^2_\beta)-tr(H_\alpha
    H_\beta)^2],\end{eqnarray*}
    and $$\sum_{i,j,k,
\alpha,\beta\neq n+1} h^\alpha_{ij}h^\beta_{ki}R_{\alpha\beta jk}=\sum_{\alpha,\beta\neq n+1} [tr(H^2_\alpha H^2_\beta)-tr(H_\alpha
    H_\beta)^2] .$$
    Since $\alpha,\beta \neq
 n+1$, $(tr(H_\alpha H_\beta))$ is a symmetric
$(p-1)\times (p-1)$-matrix. We choose the normal vector fields
$\{e_{\alpha}\}_{\alpha\neq n+1}$ such that
$$tr(H_\alpha H_\beta)= trH^2_\alpha \cdot \delta_{\alpha\beta},$$
 which implies
\begin{equation} \sum_{\alpha,\beta\neq n+1} [tr(H_\alpha H_\beta)]^2=
 \sum_{\alpha\neq n+1}tr(H^2_\alpha)^2.\end{equation}
For any real number $a$, we have
\begin{eqnarray} \sum_{i,j, \alpha\neq n+1} h^\alpha_{ij}\Delta
h^\alpha_{ij}&=&(1+a)\sum_{i,j,k,m,\alpha\neq n+1}(
h^\alpha_{ij}h^\alpha_{km}R_{mijk}+
h^\alpha_{ij}h^\alpha_{mi}R_{mkjk})-ancS_I\nonumber\nonumber\\
    &&+(a-1)\sum_{\alpha,\beta\neq n+1} [tr(H^2_\alpha H^2_\beta)-tr(H_\alpha H_\beta)^2]
     + a\sum_{\alpha\neq n+1}(trH_\alpha^2)^2\nonumber\\
     &&+a\Big\{-\sum_{\alpha\neq n+1} tr(H^2_\alpha H_{n+1})\cdot trH_{n+1}
    +\sum_{\alpha\neq n+1} [tr(H_\alpha H_{n+1})]^2\Big\}.
\end{eqnarray}
\hspace*{5mm}When $p=1$, $M$ is a compact hypersurface with nonzero
constant mean curvature and nonnegative sectional curvature in
$F^{n+1}(c)$.
  The assertion was proved by Nomizu and Symth \cite{Nomizu} for $c\geq 0$ and by Walter \cite{Walter} for $c<0$, respectively.\\
\hspace*{5mm}When $p=2$, $K_{M}\geq 0$ and $H=constant\neq0$. We
know from Theorem 9 in \cite{Yau} that $M$ is a minimal hypersurface
in the totally umbilical sphere
$S^{n+1}\Big(\frac{1}{\sqrt{c+H^2}}\Big)$. This together with
Theorem A implies that $M$ is either a totally umbilical sphere or
the standard immersion of the product of
two spheres. \\
\hspace*{5mm}When $p\geq 3$, it follows from Lemma 1 and  the
assumption that $M$ is pseudo-umbilical, i.e., $ h^{n+1} _{ij}= H
 \delta_{ij}$. Hence, we have
\begin{eqnarray}
 &&\sum_{\alpha\neq n+1} tr(H^2_\alpha H_{n+1})\cdot trH_{n+1}
    -\sum_{\alpha\neq n+1} [tr(H_\alpha H_{n+1})]^2\nonumber \\
&=& \sum_{i,j,k,m,\alpha\neq n+1}h^\alpha _{ij} h^\alpha
_{mi}h^{n+1} _{mj}h^{n+1} _{kk}-\sum_{i,j,k,m,\alpha\neq
n+1} h^\alpha _{ij}h^\alpha _{km}
    h^{n+1} _{mk}h^{n+1} _{ij}\nonumber\\
&=& nH^2 \sum_{i,j,\alpha\neq n+1}(h^\alpha _{ij})^2-H^2\sum_{\alpha\neq n+1} (trH_\alpha)^2\nonumber\\
&=& nH^2S_I.
\end{eqnarray}
On the other hand, we get from (12)
\begin{equation}\sum_{i,j,k,m,\alpha \neq n+1}h^\alpha_{ij}h^\alpha_{km}R_{mijk}+
    \sum_{i,j,k,m,\alpha \neq n+1}
    h^\alpha_{ij}h^\alpha_{mi}R_{mkjk} \geq nK_{\min}S_I.\end{equation}
    By a direct computation and the DDVV inequality, we obtain
    \begin{eqnarray}\sum_{\alpha,\beta\neq n+1}tr(H^2_\alpha H^2_\beta)
    -tr(H_\alpha
    H_\beta)^2&=&\frac{1}{2}\sum_{\alpha,\beta\neq n+1}tr(H_{\alpha} H_{\beta}-H_{\beta}
H_{\alpha})^{2}\nonumber\\
&\leq&\frac{1}{2}\Big(\sum_{\alpha\neq n+1}trH_{\alpha}^{2}\Big)^2\nonumber\\
&=&\frac{1}{2} S_{I}^2.\end{eqnarray}
It follows from (19), (20), (21) and (22) that
 \begin{eqnarray}
 \frac{1}{2}\Delta S_{I}&=&\sum_{i,j,\alpha\neq n+1} (h^\alpha_{ijk})^2 +\sum_{i,j,\alpha\neq n+1} h^\alpha_{ij}\Delta
    h^\alpha_{ij}\nonumber\\
 &\geq&(1+a)nK_{\min}S_I+a\sum_{\alpha\neq n+1}(trH_\alpha ^2)^2+
  \frac{1}{2}(a-1)S_I^2-an(c+H^2)S_I \nonumber\\
 &\geq
 &(1+a)nK_{\min}S_I+\Big(\frac{a}{p-1}+\frac{a-1}{2}\Big)S_I^2-an(c+H^2)S_I\nonumber\\
 &= & S_I\Big[(1+a)nK_{\min}+\Big(\frac{a}{p-1}+\frac{a-1}{2}\Big)S_I-an(c+H^2)\Big],
\end{eqnarray}
for $0\leq a<1$. Taking $a=\frac{p-1}{p+1}$, we get
 \begin{eqnarray*}
   \frac{1}{2}\Delta S_{I}
&\geq &nS_I[(1+a)K_{\min}-a(c+H^2)]\\
&=&nS_I\Big[\Big(1+\frac{p-1}{p+1}\Big)K_{\min}-\frac{p-1}{p+1}(c+H^2)\Big].
\end{eqnarray*}
It follows from the assumption and the maximum principal that
$S_{I}$ is a
 constant, and
$$S_I\Big[\Big(1+\frac{p-1}{p+1}\Big)K_{\min}-\frac{p-1}{p+1}(c+H^2)\Big]=0.$$

If there is a point $q\in M$ such that
 $K_{\min}(q)>\frac{(p-1)(c+H^2)}{2p},$ then
$S_{I}=0$, i.e., $M$ is a compact hypersurface with nonzero constant
mean curvature and positive sectional curvature in a totally geodesic submanifold $F^{n+1}(c)$. Therefore, $M$ is a totally umbilical sphere
$S^{n}(\frac{1}{\sqrt{c+H^2}})$.

If $K_{\min}\equiv\frac{(p-1)(c+H^2)}{2p},$ then inequalities in
(21), (22) and (23) become equalities. This together with the DDVV
inequality implies that $p=3$ and $K_{\min}=\frac{c+H^2}{3}$. Taking
$a=0$ in $(23)$, we get $S_{I}=\frac{2n}{3}(c+H^2)$. By the same
argument as in \cite{Chern}, we conclude that $n=2$. Hence,
$K_{M}=\frac{c+H^2}{3}$ and $M$ is the Veronese surface
in  $S^4(\frac{1}{\sqrt{c+H^2}})$. This  completes the proof of Theorem 2.\\\\
\hspace*{5mm} Combing Theorems 1, 2 and rigidity results in
\cite{Itoh2, Shen, Xu2}, we present a general version of the Yau
rigidity
theorem.\\\\
\textbf{Generalized Yau Rigidity Theorem.}\emph{ Let $M^n$ be an
$n$-dimensional oriented
compact submanifold with parallel mean curvature in $F^{n+p}(c)$, where $c+H^2>0$. Set $k(m,n)=\min\{sgn(m-1)m, n\}.$ Then we have\\\\
$(i)$ if $H=0$ and
$$K_{M}\geq\frac{k(p,n)c}{2[k(p,n)+1]},$$ then M is either a
totally geodesic sphere, the standard immersion
of the product of two spheres, or the Veronese submanifold;\\\\
$(ii)$ if $H \neq 0$ and $$K_M \ge
\frac{k(p-1,n)(c+H^2)}{2[k(p-1,n)+1]},$$ then $M$ is either a
totally umbilical sphere $S^n(\frac{1}{\sqrt{c+H^2}})$ in
$F^{n+p}(c)$, the standard  immersion of the product of two spheres,
or the Veronese submanifold.}\\\\
\hspace*{5mm} Recently Andrews and Baker \cite{Andrews} generalized
a weaker version of Huisken's convergence theorem \cite{Huisken} for
mean curvature flow of convex hypersurfaces in $\mathbf{R}^{n+1}$ to
higher codimensional cases. Motivated by Generalized Yau Rigidity
Theorem, we would like to propose the following conjecture on mean
curvature flow in higher codimensions,
which can be considered as a generalization of the Huisken convergence theorem \cite{Huisken}.\\\\
\textbf{Conjecture.} \emph{Let $M_0=F_0(M)$ be an $n$-dimensional
compact submanifold in an $(n+p)$-dimensional space form
$F^{n+p}(c)$ with $c+H^2>0$. If the sectional curvature of $M_0$
satisfies
$$K_M >
\frac{k(p,n)(c+H^2)}{2[k(p,n)+1]},$$ then the mean curvature flow
\begin{eqnarray}
\label{MCF}\left\{
\begin{array}{ll}
\frac{\partial}{\partial t}F(x,t)=n\xi(x,t), \,\, x\in M, \, t\ge0,  \\
F(\cdot,0)=F_0(\cdot).
\end{array}\right.
\end{eqnarray}
has a unique smooth solution $F : M \times [0, T) \rightarrow
F^{n+p}(c)$ on a finite maximal time interval, and $F_t(\cdot)$
converges uniformly to a round point
$q\in F^{n+p}(c)$ as $t \rightarrow T$.}\\\\
\hspace*{5mm}When $p=1$ and $c=0$, the conjecture was verified by
Huisken \cite{Huisken}. When $p=1$ and $c=1$, a weaker version of
the conjecture was proved by Huisken \cite{Huisken2}.

Juan-Ru Gu

Center of Mathematical Sciences\

Zhejiang University\

Hangzhou 310027\

China

E-mail address: gujr@cms.zju.edu.cn\\\\

Hong-Wei Xu

Center of Mathematical Sciences\

Zhejiang University\

Hangzhou 310027\

China

E-mail address: xuhw@cms.zju.edu.cn


\begin{thebibliography}{bb}

\bibitem{Andrews} B. Andrews and C. Baker, Mean curvature flow of pinched submanifolds
to spheres, \emph{J. Differential Geom.}, {\bf85}(2010), 357-396.
\bibitem{Chern}S. S. Chern, M. do Carmo and S. Kobayashi, Minimal submanifolds of a sphere with second fundamental form of constant length, in Functional Analysis and Related Fields, Springer-Verlag, New York(1970).
\bibitem{Choi}T. Choi and Z. Lu, On the DDVV conjecture and the comass in calibrated geometry (I), \emph{Math. Z.}, {\bf260}(2008), 409-429.
\bibitem{DeSmet}P. J. De Smet, F. Dillen, L. Verstraelen and L. Vrancken, A pointwise inequality in
submanifold theory, \emph{Arch. Math.}, {\bf 35}(1999), 115-128.
\bibitem{Ejiri}N. Ejiri, Compact minimal submanifolds of a sphere with positive Ricci curvature, \emph{J. Math.
Soc. Japan.}, {\bf31}(1979), 251-256.
\bibitem{Erbacher}J. Erbacher, Reduction of the codimension of an isometric
immersion. \emph{J. Differential Geom.}, {\bf5}(1971), 333-340.
\bibitem{Ge} J. Q. Ge and Z. Z. Tang, A proof of the DDVV conjecture and its equality case, \emph{Pacific J. Math.},
 {\bf237}(2008), 87-95.
\bibitem{Huisken} G. Huisken, Flow by mean curvature of convex surfaces into
spheres,  \emph{J. Differential Geom.}, {\bf20}(1984), 237¨C266.
\bibitem{Huisken2} G. Huisken, Deforming hypersurfaces of the sphere by their mean
curvature, \emph{Math. Z.}, {\bf195}(1987), 205¨C219.
\bibitem{Kozlowski} M. Kozlowski and U. Simon, Minimal immersions of $2$-manifolds into spheres, \emph{Math. Z.}, {\bf186}(1984), 377--382.
\bibitem{Itoh1}T. Itoh, On veronese manifolds, \emph{J. Math. Soc. Japan.}, {\bf 27}(1975), 497-506.
\bibitem{Itoh2}T. Itoh, Addendum to my paper "On veronese manifolds", \emph{J. Math. Soc. Japan.}, {\bf 30}(1978), 73-74.
\bibitem{Lawson}B. Lawson, Local rigidity theorems for minimal hyperfaces, \emph{Ann. of Math.}, {\bf 89}(1969), 187-197.
\bibitem{Li}A. M. Li and J. M. Li, An intrinsic rigidity theorem for minimal submanifolds in a sphere, \emph{Arch. Math.}, {\bf 58}(1992), 582-594.
\bibitem {Lu1}Z. Lu, On the DDVV conjecture and the comass in calibrated geometry (II), arXiv:math.DG/0708.2921v1.
\bibitem {Lu2}Z. Lu, Proof of the normal scalar curvature conjecture, arXiv:math.DG/0711.3510v1.
\bibitem {Lu3}Z. Lu, Recent developments of the DDVV conjecture, \emph{Bull. Transil. Univ. Brasov serB.}, {\bf14}(2008), 133-144.
\bibitem{Lu4}Z. Lu, Normal scalar curvature conjecture and its applications, arXiv:math.DG/ 0803.0502.
\bibitem{Nomizu}K. Nomizu and B. Smyth, A formula of Simons' type and hypersurfaces with constant mean curvature, \emph{J. Differential Geom.}, {\bf 3}(1969), 367-377.
\bibitem{Shiohama}K. Shiohama and H. W. Xu, A general rigidity theorem for complete submanifolds, \emph{Nagoya Math. J.}, {\bf 150}(1998), 105-134.
\bibitem{Shen}Y. B. Shen, Submanifolds with nonnegative sectional curvature, \emph{Chinese Ann. Math. SerB}, {\bf 5}(1984), 625-632.
\bibitem{Simons}J. Simons, Minimal varieties in Riemannian manifolds, \emph{Ann. Math.}, {\bf 88}(1986), 62-105.
\bibitem{Walter}R. Walter, Compact hypersurfaces with a constant higher mean curvature function, \emph{Math. Ann.,} {\bf270}(1985), 125-145.
\bibitem{Xu}H. W. Xu, A rigidity theorem for submanifolds with parallel mean curvature in a sphere, \emph{Arch. Math.}, {\bf 61}(1993), 489-496.
\bibitem{Xu1}H. W. Xu, On closed minimal submanifolds in pinched Riemannian
manifolds, \emph{Trans. Amer.  Math. Soc.}, {\bf 347}(1995),
1743-1751.
\bibitem{Xu2}H. W. Xu and W. Han, Geometric rigidity theorem for submanifolds with positive curvature, \emph{Appl. Math. J. Chinese Univ. Ser. B},
 {\bf 20}(2005), 475-482.
\bibitem{Yau}S. T. Yau, Submanifolds with constant mean curvature I, II, \emph{Amer. J. Math.}, {\bf96, 97}(1974, 1975), 346-366, 76-100.

\end{thebibliography}
\end{document}